\documentclass{article}
\usepackage{amsmath,amssymb,amsthm}
\begin{document}

\title{ Fourier-Borel transformation \\ on the hypersurface of any reduced polynomial}
\author{ Atsutaka \textsc{Kowata} and Masayasu \textsc{Moriwaki} }
\date{}
\maketitle

\begin{abstract}
For any polynomial $p$ on $\mathbf{C}^{n}$, 
a variety $ V_{p} = \{ z \in \mathbf{C}^{n} ; p(z)=0 \} $ 
will be considered. 
Let $\text{Exp}(V_{p})$ be the space of holomorphic functions of 
expotential growth on $V_{p}$. 
We shall prove that 
the Fourier-Borel transformation yields an isomorphism 
of the dual space $\text{Exp}'(V_{p})$ with the space of holomorphic solutions 
$\mathcal{O}_{\partial p}(\mathbf{C}^{n})$ 
with respect to the differential operator $\partial p$ which is obtained by replacing 
each variable $z_{j}$ with $\partial / \partial z_{j}$ in $p$ 
when $p$ is a reduced polynomial. 
The result has been shown by Morimoto  
and by Morimoto-Wada-Fujita only for the case 
$p(z) = z_{1}^{2} + \dots + z_{n}^{2} + \lambda \, (n \geq 2)$. 


\end{abstract}

\textbf{\S 1. Introduction and Preliminaries.}

Let $\mathcal{O}(\mathbf{C}^{n})$ be the space of entire functions on $\mathbf{C}^{n}$ 
equipped with the topology of uniform convergence on compact subsets. 
$\mathcal{O}(\mathbf{C}^{n})$ is a FS (Fr\'{e}chet-Schwartz) space. 
We put
\begin{gather*}
   \| f \|_{A} = \sup \{ |f(z)| \exp (-A|z|) \ ; \ z \in \mathbf{C}^{n} \} 
   \ \ \text{and} \\ 
   E_{A} = \{ f \in \mathcal{O}(\mathbf{C}^{n}) \ ; \ \| f \|_{A} < \infty \} 
\end{gather*}
for any positive number $A$, 
where $|z|= \sqrt{|z_{1}|^{2} + \dots + |z_{n}|^{2} }$ for $z=(z_{1},\dots,z_{n}) \in \mathbf{C}^{n}$. 
The space $E_{A}$ is a Banach space with respect to the norm $\| \ \|_{A}$. 
We define a topological linear space 
$ \text{Exp}(\mathbf{C}^{n}) = \text{ind lim}_{A>0} \ E_{A}  $ 
equipped with the inductive limit topology. 
As is well known, $\text{Exp}(\mathbf{C}^{n})$ is a DFS space 
(a dual Fr\'{e}chet-Schwartz space) 
and called the space of entire functions of exponential type. 
We denote the dual space of $\text{Exp}(\mathbf{C}^{n})$ by $\text{Exp}'(\mathbf{C}^{n})$. 
It is clear that $\text{Exp}'(\mathbf{C}^{n})$ becomes a FS space by the srtong dual topology. 

Moreover it is easily seen that 
$f+g$, $fg \in \text{Exp}(\mathbf{C}^{n})$ for any 
$f, \, g \in \text{Exp}(\mathbf{C}^{n})$, that is, 
$\text{Exp}(\mathbf{C}^{n})$ is a commutative algebra 
with respect to the usual sum and product of functions. 

For any $f,g \in \text{Exp}(\mathbf{C}^{n})$ and $T \in \text{Exp}'(\mathbf{C}^{n})$, 
we define $(gT)(f)=T(gf)$. 
Since $\text{Exp}(\mathbf{C}^{n})$ is a commutative algebra and $gT$ is 
a continuous linear functional, 
we have $gT \in \text{Exp}'(\mathbf{C}^{n})$.
\\

\textsc{Definition 1.1.} \ For any $T \in \text{Exp}'(\mathbf{C}^{n})$, 
we define the Fourier-Borel transformation $\mathcal{F}$ by 
\[ \mathcal{F}(T)(z) = \left\langle T, \, e^{\langle z, \, \cdot \, \rangle} 
\right\rangle, \]
where $\langle z,w \rangle = z_{1} w_{1} + \dots + z_{n} w_{n}$ 
for $z,w \in \mathbf{C}^{n}$ 
and $\langle T , f \rangle$ is the dual pairing: 
$\langle T , f \rangle = T(f)$ for any $T \in \text{Exp}'(\mathbf{C}^{n})$ 
and $f \in \text{Exp}(\mathbf{C}^{n})$. 
\\

For any polynomial $p$ on $\mathbf{C}^{n}$, 
we set a variety $ V_{p} = \{ z \in \mathbf{C}^{n} ; p(z)=0 \} $. 
$ V_{p} $ is a closed set of $\mathbf{C}^{n} .$  
Thanks to Oka-Cartan's Theorem, 
the restriction mapping $ r : \mathcal{O}(\mathbf{C}^{n}) \to \mathcal{O}(V_{p})$ is surjective. 
Hence, we have the following exact sequence:
   \[ 0 \to \mathcal{K}_{p} \stackrel{i}{\longrightarrow} \mathcal{O}(\mathbf{C}^{n})
      \stackrel{r}{\longrightarrow} \mathcal{O}(V_{p}) \to 0 , \]
where $ \mathcal{K}_{p} = \{ f \in \mathcal{O}(\mathbf{C}^{n}) ; f|_{V_{p}}=0 \} $ 
is a closed subspace of $\mathcal{O}(\mathbf{C}^{n})$ and 
$i$ is the canonical injection. 

We define the space $\text{Exp}(V_{p})$ by the image of the space $\text{Exp}(\mathbf{C}^{n})$ 
of the entire functions of exponential type under the restriction mapping $r$. 
The topology of $\text{Exp}(V_{p})$ is defined by the quotient topology of the restriction mapping $r$. 
We set $ \mathcal{K}_{p}^{E} = \mathcal{K}_{p} \cap \text{Exp}(\mathbf{C}^{n}) $. 
$\mathcal{K}_{p}^{E}$ is a closed subspace of $\text{Exp}(\mathbf{C}^{n})$. 
By definition, we have the exact sequence 
   \[ 0 \to \mathcal{K}_{p}^{E} \stackrel{i}{\longrightarrow} \text{Exp}(\mathbf{C}^{n})
      \stackrel{r}{\longrightarrow} \text{Exp}(V_{p}) \to 0 \]
and $ \text{Exp}(V_{p}) \cong \text{Exp}(\mathbf{C}^{n}) / \mathcal{K}_{p}^{E} .$ 
Hence $\text{Exp}(V_{p})$ is a DFS space, being a quotient space of a DFS space 
by a closed subspace. 

Let $\text{Exp}'(V_{p})$ be the dual space of $\text{Exp}(V_{p})$. 
The space $\text{Exp}'(V_{p})$ becomes a FS space by the strong dual topology, 
since $\text{Exp}(V_{p})$ is a DFS space. 
Because the restriction mapping $ r:\text{Exp}(\mathbf{C}^{n}) \to \text{Exp}(V_{p}) $ is surjective, 
the transposed mapping $ {}^{t}r: \text{Exp}'(V_{p}) \to \text{Exp}'(\mathbf{C}^{n}) $ is injective. 

Let $\partial p$ be a differential operator obtained 
by replacing each variable $z_{j}$ 
with $\partial / \partial z_{j}$ in $p$. 
We set a holomorphic solution space 
$ \mathcal{O}_{\partial p}(\mathbf{C}^{n}) = \{ f \in \mathcal{O}(\mathbf{C}^{n}) ; 
   \partial p(f)=0 \} $. 
Since the mapping $\partial p : \mathcal{O}(\mathbf{C}^{n}) \to \mathcal{O}(\mathbf{C}^{n})$ 
is continuous, $\mathcal{O}_{\partial p}(\mathbf{C}^{n})$ 
is a closed subspace of the FS space $\mathcal{O}(\mathbf{C}^{n})$. 
Thus $\mathcal{O}_{\partial p}(\mathbf{C}^{n})$ is a FS space. 

The purpose of this paper is to give the following theorem; 
\\

\textsc{Theorem 1.2.} \  
\emph{The composed mapping} 
\[ \mathcal{F} \circ {}^{t}r : \text{Exp}'(V_{p}) 
\to \mathcal{O}_{\partial p}(\mathbf{C}^{n}) \]
\emph{is a topological linear isomorphism, 
if and only if $p$ is a reduced polynomial on $\mathbf{C}^{n}$. 
We will abbreviate $\mathcal{F} \circ {}^{t}r$ to $\mathcal{F}$.}
\\

Here we recall the definition of \emph{reduced polynomial}. 
If a principal ideal $\langle p \rangle$ in a polynomial ring on $\mathbf{C}^{n}$ 
generated by a polynomial $p$ 
is a reduced ideal, 
that is, $\langle p \rangle = \sqrt{\langle p \rangle}$, 
$p$ is called a reduced polynomial. 
A reduced polynomial is represented by 
the product of irreducible polynomials which has no multiplicity. 
An irreducible polynomial is obviously a reduced polynomial. 

We show the above theorem by means of Oka-Cartan's and 
Martineau's theorems. 
Before giving the proof, we explain some known results. 

For a polynomial $p(z) = z_{1}^{2} + \dots + z_{n}^{2} + \lambda \ 
(n \geq 2, \, \lambda \neq 0)$, 
we see that 
$\partial p = \varDelta_{z}+ \lambda$, where 
$\varDelta_{z} = \partial^{2} / \partial z_{1}^{2} + \dots + 
 \partial^{2} / \partial z_{n}^{2}$ is called the complex Laplacian, 
and 
$V_{p}$ is isomorphic to the complex sphere $\tilde{S}^{n-1}$ 
defined by $\{ z \in \mathbf{C}^{n} ; z_{1}^{2} + \dots + z_{n}^{2} = 1 \}$. 
Since $p$ is an irreducible polynomial, 
Theorem 1.2 implies 
\\

\textsc{Theorem 1.3} (Morimoto \cite{mo83} \cite{mo93} \cite{mo98}). \ 
\emph{The Fourier-Borel transformation}  
   \begin{equation*}
      \mathcal{F} : \text{Exp}'(\tilde{S}^{n-1}) \stackrel{\sim}{\longrightarrow} 
      \mathcal{O}_{\lambda} (\mathbf{C}^{n})
   \end{equation*}
\emph{is a topological linear isomorphism, 
where $\mathcal{O}_{\lambda}(\mathbf{C}^{n})$ is the space of eigenfunctions} 
$\{ f \in \mathcal{O}(\mathbf{C}^{n}) ; (\varDelta_{z}+ \lambda)f = 0 \}$ 
\emph{with respect to the eigenvalue} $- \lambda$. 
\\

For a polynomial $p(z) = z_{1}^{2} + \dots + z_{n}^{2} \, (n \geq 2)$, 
we see that $\partial p = \varDelta_{z}$ and 
$V_{p} = \{ z \in \mathbf{C}^{n} ; z_{1}^{2}+\dots+z_{n}^{2} =0 \}$, 
where $V_{p}$ is called the complex light cone. 
$p$ is an irreducible polynomial for $n \geq 3$ 
and still a reduced polynomial for $n=2$. 
Hence, Theorem 1.2 implies
\\

\textsc{Theorem 1.4} 
(Morimoto-Wada-Fujita \cite{mo93} \cite{mo98} \cite{m-f} \cite{m-w}). \ 
\emph{The Fourier-Borel transformation} 
   \begin{equation*}
      \mathcal{F} : \text{Exp}'(\text{complex light cone}) 
      \stackrel{\sim}{\longrightarrow} 
      \mathcal{O}_{\varDelta_{z}} (\mathbf{C}^{n})
   \end{equation*}
\emph{is a topological linear isomorphism, 
where $\mathcal{O}_{\varDelta_{z}}(\mathbf{C}^{n})$ is the space} 
$\{ f \in \mathcal{O}(\mathbf{C}^{n}) ; \varDelta_{z} f = 0 \} $. 
\\

\textbf{\S 2. Isomorphism given by Fourier-Borel transformation.}

The following theorem plays an important role in this section. 
\\

\textsc{Theorem 2.1} (Martineau~\cite{ma}). \  
\emph{The Fourier-Borel transformation}
\[ \mathcal{F} : \text{Exp}'(\mathbf{C}^{n}) 
\longrightarrow \mathcal{O}(\mathbf{C}^{n}) \]
\emph{is a topological linear isomorphism.} 
\\

It is clear that for any $g \in \text{Exp}(\mathbf{C}^{n}) $, the mapping 
\[ \tau_{g} : \text{Exp}'(\mathbf{C}^{n}) \ni 
T \mapsto g T \in \text{Exp}'(\mathbf{C}^{n}) \]
is linear and continuous. 
We set a subspace 
\[ \text{Exp}'(\mathbf{C}^{n})_{p} 
   = \{ T \in \text{Exp}'(\mathbf{C}^{n}) ; \tau_{p}(T)=pT=0 \} , \]
that is, $\text{Exp}'(\mathbf{C}^{n})_{p} = \ker \tau_{p}$. 
$\text{Exp}'(\mathbf{C}^{n})_{p}$ 
is a FS space as a closed subspace of the FS space $\text{Exp}'(\mathbf{C}^{n})$. 

Owing to Martineau's theorem, we have the following proposition. 
\\

\textsc{Proposition 2.2.} \ 
\emph{The following restriction, 
denoted by the same notation $\mathcal{F}$, 
of the Fourier-Borel transformation above}:
   \begin{equation*}
      \mathcal{F} : \text{Exp}'(\mathbf{C}^{n})_{p} \stackrel{\sim}{\longrightarrow} 
      \mathcal{O}_{\partial p} (\mathbf{C}^{n})
   \end{equation*}
\emph{is a topological linear isomorphism.} 

\begin{proof}[\textsc{Proof.}]
Obviously, the mappings
$ \tau_{p} : \text{Exp}'(\mathbf{C}^{n}) \ni T 
\mapsto pT \in \text{Exp}'(\mathbf{C}^{n}) $ 
and $ \partial p : \mathcal{O}(\mathbf{C}^{n}) \to \mathcal{O}(\mathbf{C}^{n}) $ 
are continuous. 
Moreover it is easily seen that the following diagram commutes: 
\[ \begin{array}{ccc}
      \text{Exp}'(\mathbf{C}^{n}) & 
      \begin{matrix} \stackrel{\sim}{\longrightarrow} \\  \mathcal{F} \end{matrix} & 
      \mathcal{O}(\mathbf{C}^{n}) \ \ \ \ \ \\
      \tau_{p} \ \Big\downarrow  & \circlearrowleft & \Big\downarrow \ \partial p \ \ \ \\
      \text{Exp}'(\mathbf{C}^{n}) & 
      \begin{matrix} \stackrel{\sim}{\longrightarrow} \\  \mathcal{F} \end{matrix} & 
      \mathcal{O}(\mathbf{C}^{n}) .  \ \ \ \ \
   \end{array} \]
by the Fourier-Borel transformation $\mathcal{F}$. 
So, $\ker \tau_{p} \simeq \ker \partial p$. 
\end{proof}

We set a subspace 
   \[ \text{Exp}'(\mathbf{C}^{n};\mathcal{K}_{p}^{E}) =
      \{ T \in \text{Exp}'(\mathbf{C}^{n});T|_{\mathcal{K}_{p}^{E}}=0 \} , \]
where the mapping $T|_{\mathcal{K}_{p}^{E}}$ is 
the restriction of the linear mapping $T$ over the subspace $\mathcal{K}_{p}^{E}$. 
It is obvious that $\text{Exp}'(\mathbf{C}^{n};\mathcal{K}_{p}^{E})$ is a closed subspace of 
$\text{Exp}'(\mathbf{C}^{n}) $. Indeed, let $i : \mathcal{K}_{p}^{E} \to \text{Exp}(\mathbf{C}^{n})$
be the canonical injection, then we have $T|_{\mathcal{K}_{p}^{E}}= {}^{t}i(T)$. 
Thus we have $ \text{Exp}'(\mathbf{C}^{n};\mathcal{K}_{p}^{E}) = \ker {}^{t}i $. 
Therefore $\text{Exp}'(\mathbf{C}^{n};\mathcal{K}_{p}^{E})$ becomes a FS space. 
\\

\textsc{Proposition 2.3.} \ 

(1). \emph{The transposed mapping} $ {}^{t}r: \text{Exp}'(V_{p}) \to 
           \text{Exp}'(\mathbf{C}^{n};\mathcal{K}_{p}^{E}) $ 
     \emph{is a topological linear isomorphism and} 
           $ \text{Exp}'(\mathbf{C}^{n};\mathcal{K}_{p}^{E})$ 
     \emph{is a subspace of} 
           $\text{Exp}'(\mathbf{C}^{n})_{p} $. 

(2). \emph{If $\mathcal{K}_{p}^{E}$ is a principal ideal of} 
           $\text{Exp}(\mathbf{C}^{n})$ 
     \emph{generated by $p$, then we have}
            \[ \text{Exp}'(\mathbf{C}^{n};\mathcal{K}_{p}^{E})=
            \text{Exp}'(\mathbf{C}^{n})_{p}. \]

\begin{proof}[\textsc{Proof.}] 
(1). It is easily seen that the transposed mapping ${}^{t}r$ is linear, continuous and injective. 
Indeed, for any $S \in \text{Exp}'(V_{p})$ and $f \in \mathcal{K}_{p}^{E}$, we have 
   \[ \langle {}^{t}r(S),f \rangle = \langle S,r(f) \rangle = 0. \]
This implies that  
   \[ {}^{t}r(\text{Exp}'(V_{p})) \subset \text{Exp}'(\mathbf{C}^{n};\mathcal{K}_{p}^{E}) .\]
Let $T$ be an element of $\text{Exp}'(\mathbf{C}^{n};\mathcal{K}_{p}^{E})$. 
Since $r : \text{Exp}(\mathbf{C}^{n}) \to \text{Exp}(V_{p})$ 
is surjective and $\ker r \subset \ker T$, 
there exists unique linear mapping 
$ S:\text{Exp}(V_{p}) \to \mathbf{C} $ such that $T=S \circ r$. 
If $U$ is an open subset of $\mathbf{C}$, then $r(T^{-1}(U))=S^{-1}(U)$ 
since $r$ is surjective. On the other hand, because $r$ is an open mapping, 
$S$ is a continuous mapping. Hence the mapping $S$ belongs to $\text{Exp}'(V_{p})$. 
Moreover, since ${}^{t}r(S)=T$, we obtain the surjectivity of ${}^{t}r$. 
By the closed graph theorem for FS spaces, we get the first assertion.
The second assertion is clear from the definitions of 
$\text{Exp}'(\mathbf{C}^{n})_{p}$ and 
$\mathcal{K}_{p}^{E}$. 

(2). If $\mathcal{K}_{p}^{E}$ is a principal ideal of $\text{Exp}(\mathbf{C}^{n})$ generated by the polynomial $p$, 
then for $f \in \mathcal{K}_{p}^{E}$ there exists a function $ g \in \text{Exp}(\mathbf{C}^{n}) $ such that  
$f=pg$. 
So, if $T \in \text{Exp}'(\mathbf{C}^{n})_{p}$ and $f \in \mathcal{K}_{p}^{E}$
then $T(f)=T(pg)=pT(g)=0$ 
and hence $T \in \text{Exp}'(\mathbf{C}^{n};\mathcal{K}_{p}^{E})$. 
\end{proof}

From Propositions 2.2 and 2.3, we have the following corollary. 
\\

\textsc{Corollary 2.4.} \ 
\emph{Let $p$ be a polynomial on $\mathbf{C}^{n}$. 
If $\mathcal{K}_{p}^{E}$ is a principal ideal of}
$\text{Exp}(\mathbf{C}^{n})$ 
\emph{generated by $p$, then the composed mapping} 
\[ \mathcal{F} \circ {}^{t}r : \text{Exp}'(V_{p}) \to \mathcal{O}_{\partial p}(\mathbf{C}^{n}) \]
\emph{is a topological linear isomorphism. 
We will abbreviate $\mathcal{F} \circ {}^{t}r$ to $\mathcal{F}$.} 
\\

\textbf{\S 3. Proof of Theorem 1.2.}

In this section, we shall prove Theorem 1.2. 
We need some lemmas and propositions for the proof. 

First of all, we think of exponential growth in one variable case. 
Let $p$ be a polynomial defined by 
$p(z) = a_{0} z^{d} + a_{1} z^{d-1} + \dots +a_{d}$ on $\mathbf{C}$. 
We fix a complex number $\xi$ and a positive number $r$. 
Owing to the P\'{o}lya-Szeg\"{o}'s result~\cite[p.~86~problem~66]{p-s}, 
we can find some positive number $ 0 < \rho \leq r$ 
such that 
\[ 2 |a_{0}| \left( \frac{r}{4} \right)^{d} \leq |p(\zeta)|, \ \ \ 
\text{for any} \ \zeta \in \{ z \in \mathbf{C} ; | z - \xi | = \rho \} . \]
Let $f$ be a holomorphic function on $\{ z \in \mathbf{C} ; |z - \xi | \leq r \}$ 
satisfying 
\[ |p(z) f(z)| \leq M e^{A|z|}, \ \ \ \text{for some} \ A>0, \, M \geq 0. \]
Applying the maximal principle of holomorphic functions to $f$, 
there exists $\zeta_{0} \in \{ z \in \mathbf{C} ; |z-\xi| = \rho \}$ such that 
\[ |p(\zeta_{0}) f(\xi)| \leq |p(\zeta_{0}) f(\zeta_{0})|. \] 
Thus we have the following Lemma, 
putting a positive constant $c = 4^{d} / 2 r^{d}$. 
\\

\textsc{Lemma 3.1.} \ 
\emph{Fix a polynomial $p$ on $\mathbf{C}$, 
an element $\xi \in \mathbf{C}$ and an $r>0$. 
Suppose $f$ is a holomorphic function on $\{ z \in \mathbf{C} ; |z-\xi| \leq r \}$ 
satisfying $|p(z) f(z)| \leq M e^{A|z|}$ for some $A>0$ and $M \geq 0$. 
Then we have $|a_{0} f(\xi)| \leq c e^{A(|\xi|+r)} M$, 
where $c$ is a positive constant depending only on $p$ and $r$. }
\\

Next, we shall extend Lemma 3.1 to the n-variable case. 
Let $p$ be a non-zero polynomial on $\mathbf{C}^{n}$ 
and fix any $\xi = (\xi_{1}, \dots, \xi_{n})$ in $\mathbf{C}^{n}$ and any $r>0$. 
Suppose $f$ is a holomorphic function on the polydisk 
$\{ z=(z_{1}, \dots, z_{n}) \in \mathbf{C}^{n} ; |z_{i}-\xi_{i}| \leq r \ 
 (1 \leq i \leq n) \}$ 
satisfying $|p(z) f(z)| \leq M e^{A|z|}$ for some $A>0$ and $M \geq 0$. 
First, we fix $z'=(z_{1}, \dots, z_{n-1})$ in 
$\{ z' \in \mathbf{C}^{n-1}; |z_{i}-\xi_{i}| \leq r \ (1 \leq i \leq n-1) \}$, 
and regard $p$ as a polynomial of the single variable $z_{n}$ 
with degree $d$. 
Let $c$ be a positive constant $4^{d}/2r^{d}$. 
By Lemma 3.1, we have 
$|\tilde{p}(z') f(z', \xi_{n})| \leq c e^{A(|z'|+|\xi_{n}|+r)} M$, 
where $\tilde{p}(z')$ be the coefficient of $p$ with respect to $z_{n}^{d}$. 
Here we used the inequality 
$ |p(z', z_{n}) f(z', z_{n})| \leq M e^{A(|z'|+|z_{n}|)} $. 
By iteration, there exists a positive constant $\hat{c}$ depending only on $p$ and $r$ 
such that 
\[ |f(\xi)| \leq \hat{c} e^{A(|\xi_{1}| + \dots + |\xi_{n}| + nr)} M. \]
Applying Cauchy-Schwarz's inequality 
$ ( 1 \cdot |\xi_{1}|+ \dots + 1 \cdot |\xi_{n}|)^{2} \leq n |\xi|^{2} $, 
we obtain the following lemma. 
\\

\textsc{Lemma 3.2.} \ 
\emph{Fix a non-zero polynomial $p$ on $\mathbf{C}^{n}$, 
an element $\xi = (\xi_{1}, \dots, \xi_{n}) \in \mathbf{C}^{n}$ and an $r>0$. 
Suppose $f$ is a holomorphic function on the polydisk 
$\{ z=(z_{1}, \dots, z_{n}) \in \mathbf{C}^{n} ; |z_{i}-\xi_{i}| \leq r \ 
 (1 \leq i \leq n) \}$ 
satisfying $|p(z) f(z)| \leq M e^{A|z|}$ for some $A>0$ and $M \geq 0$. 
Then $|f(\xi)| e^{- \sqrt{n} A |\xi|} \leq \hat{c} e^{nrA} M$, 
where $\hat{c}$ is a positive constant depending only on $p$ and $r$. }
\\

Now, we recall the definition of 
$E_{A} = \{ f \in \mathcal{O}(\mathbf{C}^{n}) ; \| f \|_{A} < \infty \}$, 
where $\| f \|_{A} = \sup_{z \in \mathbf{C}^{n}} \{ |f(z)| e^{-A|z|} \}$. 
We have the following proposition about global exponential growth 
in $\mathbf{C}^{n}$ by Lemma 3.2. 
\\

\textsc{Proposition 3.3.} \ 
\emph{Fix a non-zero polynomial $p$ on $\mathbf{C}^{n}$ and an $A>0$. 
Suppose $F$ is an entire function satisfying $\| pF \|_{A} < \infty $. 
Then $\| F \|_{\sqrt{n}A} \leq c_{A} \| pF \|_{A}$,  
where $c_{A}$ is a positive constant depending only on $p$ and $A$. }

\begin{proof}[\textsc{Proof.}]
We fix an $r>0$. 
Since $| p(z) F(z) | \leq e^{A |z|} \| pF \|_{A}$ for any $z \in \mathbf{C}^{n}$, 
there exists a positive constant $\hat{c}$ depending only on $p$ and $r$ such that 
\[ |F(z)| \leq \hat{c} e^{ n r A } e^{\sqrt{n} A |z|} \| pF \|_{A} \] 
by Lemma 3.2. 
Setting a positive constant $c_{A} = \hat{c} e^{ n r A }$ 
depending only on $p$, $A$ and fixed positive constant $r$, 
we have 
\[ \| F \|_{\sqrt{n}A} = \sup_{z \in \mathbf{C}^{n}} |F(z)| e^{-\sqrt{n} A |z|} 
                       \leq c_{A} \| pF \|_{A}. \]
\end{proof}

We have the following proposition by Proposition 3.3. 
\\

\textsc{Proposition 3.4.} \ 
\emph{Let $p$ be a polynomial on $\mathbf{C}^{n}$. The continuous map} 
$ \sigma_{p} : \text{Exp}(\mathbf{C}^{n}) \ni f \mapsto pf 
\in \text{Exp}(\mathbf{C}^{n})$ 
\emph{is  a closed mapping.}

\begin{proof}[\textsc{Proof.}] 
Since the proposition is clear if $p \equiv 0$, 
we may assume that $p$ is a non-zero polynomial. 
Let $Z$ be a closed subset of $\text{Exp}(\mathbf{C}^{n})$. 
We take a sequence $\{ p f_{m} \}$ in $\sigma_{p}(Z)$ 
such that $pf_{m} \to g \ (m \to \infty)$ for some $g \in \text{Exp}(\mathbf{C}^{n})$. 
By the property of the inductive limit topology, 
there exists some $A>0$ such that $pf_{m} \to g \ (m \to \infty)$ in $E_{A}$. 
On the other hand, by the Proposition 3.3, 
we can see that $\{ f_{m} \}$ is a Cauchy sequence in a Banach space $E_{\sqrt{n}A}$, 
and hence we find a unique element $f$ in $E_{\sqrt{n}A}$ such that $f_{m} \to f$. 
In addition, since $Z$ is closed, $f \in Z \cap E_{\sqrt{n}A}$. 
Hence, $pf_{m} \to pf$ in $E_{\sqrt{n}A+1}$ and $pf = g$, 
because $\text{Exp}(\mathbf{C}^{n})$ is a Hausdoroff space. 
Thus the sequence $\{ pf_{m} \}$ is convergent in $\sigma_{p}(Z)$. 
Therefore $\sigma_{p}(Z)$ is closed. 
\end{proof}

\begin{proof}[\textsc{Proof of Theorem 1.2.}]  
Let $p$ be a reduced polynomial on $\mathbf{C}^{n}$ 
and $f$ an entire function such that $f|_{V_{p}} = 0$. 
Owing to R\"{u}ckert Nullstellensatz~\cite{a}, 
there exists an entire function $g$ such that $f=pg$. 
(There exists locally such a function near $V_{p}$ 
  by R\"{u}ckert Nullstellensatz, 
  which coinsides with the holomorphic function $f/p$ on $\mathbf{C}^{n}-V_{p}$.) 
Further, if $f \in \text{Exp}(\mathbf{C}^{n})$ 
then $g \in \text{Exp}(\mathbf{C}^{n})$ by Proposition 3.3. 
Thus, $\mathcal{K}_{p}^{E}= \langle p \rangle$ and 
$\text{Exp}'(V_{p}) \cong \mathcal{O}_{\partial p}(\mathbf{C}^{n})$ by Corollary 2.4, 
where $\langle p \rangle$ is a principal ideal generated by $p$, that is, 
a subspace $\{ fp ; f \in \text{Exp}(\mathbf{C}^{n}) \}$. 

Conversely, if $p$ is not a reduced polynomial, 
we can find some irreducible polynomial $p_{1}$ such that $p=p_{1}^{2} p_{2}$. 
Set $q=p_{1} p_{2}$. 
Obviously, $V_{p} = V_{q}$ and 
$\langle p \rangle \subsetneq \langle q \rangle \subset \mathcal{K}_{p}^{E}$. 
By Proposition 3.4, $\langle p \rangle$ and $\langle q \rangle$ are 
closed subspaces of DFS space  $\text{Exp}(\mathbf{C}^{n})$, 
and each space is a DFS space. 
We can choose a non-zero continuous linear map $S : \langle q \rangle \to \mathbf{C}$ 
such that $S|_{\langle p \rangle} = 0$. 
Indeed, for example, for $v \in V_{p_{1}}$, 
we define a linear map $T_{v} : \langle q \rangle \to \mathbf{C}$ by $T_{v}(fq)=f(v)$ 
for $f \in \text{Exp}(\mathbf{C}^{n})$. Fix any $A>0$. 
By Proposition 3.3, there exists some positive constant $c_{A}$ 
such that $|T_{v}(fq)| = |f(v)| \leq c_{A} e^{\sqrt{n} A |v|} \| fq \|_{A}$. 
This means that $T_{v}$ is a continuous map. 
It is clear that $T_{v} \neq 0$ and $T_{v}|_{\langle p \rangle}=0$. 

Applying Hahn-Banach's Theorem, we have $\hat{S} \in \text{Exp}'(\mathbf{C}^{n})$  satisfying $\hat{S}|_{\langle q \rangle} = S$. 
It is clear that $\hat{S} \in \text{Exp}'(\mathbf{C}^{n})_{p}$ 
and $\hat{S} \notin \text{Exp}'(\mathbf{C}^{n} : \mathcal{K}_{p}^{E})$. 
Thus, $\text{Exp}'(V_{p}) \not\simeq \mathcal{O}_{\partial p}(\mathbf{C}^{n})$ 
by Propositions 2.2 and 2.3. 
\end{proof}

Department of Mathematics, Graduate School of Science, 
Hiroshima University. 

Higashi-Hiroshima, Hiroshima 739-8526 Japan. 

e-mail: kowata@math.sci.hiroshima-u.ac.jp (A. Kowata), 

\hspace*{1cm} d042540@math.sci.hiroshima-u.ac.jp (M. Moriwaki)

\end{document}